\documentstyle[12pt]{article}
\def\hybrid{\topmargin 0pt      \oddsidemargin 0pt
          \headheight 0pt \headsep 0pt
          \voffset=-0.5cm
          \textwidth 6.25in       % A4 paper
          \textheight 9.5in       % A4 paper
          \marginparwidth 0.0in
          \parskip 5pt plus 1pt   \jot = 1.5ex}
\catcode`\@=11
\def\marginnote#1{}

\newcount\hour
\newcount\minute
\newtoks\amorpm
\hour=\time\divide\hour by60
\minute=\time{\multiply\hour by60 \global\advance\minute by-\hour}
\edef\standardtime{{\ifnum\hour<12 \global\amorpm={am}%
          \else\global\amorpm={pm}\advance\hour by-12 \fi
          \ifnum\hour=0 \hour=12 \fi
          \number\hour:\ifnum\minute<10 0\fi\number\minute\the\amorpm}}
\edef\militarytime{\number\hour:\ifnum\minute<10 0\fi\number\minute}

\def\draftlabel#1{{\@bsphack\if@filesw {\let\thepage\relax
     \xdef\@gtempa{\write\@auxout{\string
        \newlabel{#1}{{\@currentlabel}{\thepage}}}}}\@gtempa
     \if@nobreak \ifvmode\nobreak\fi\fi\fi\@esphack}
          \gdef\@eqnlabel{#1}}
\def\@eqnlabel{}
\def\@vacuum{}
\def\draftmarginnote#1{\marginpar{\raggedright\scriptsize\tt#1}}

\def\draftlabel#1{{\@bsphack\if@filesw {\let\thepage\relax
     \xdef\@gtempa{\write\@auxout{\string
        \newlabel{#1}{{\@currentlabel}{\thepage}}}}}\@gtempa
     \if@nobreak \ifvmode\nobreak\fi\fi\fi\@esphack}
          \gdef\@eqnlabel{#1}}
\def\@eqnlabel{}
\def\@vacuum{}
\def\draftmarginnote#1{\marginpar{\raggedright\scriptsize\tt#1}}

\def\draft{\oddsidemargin -.5truein
          \def\@oddfoot{\sl preliminary draft \hfil
          \rm\thepage\hfil\sl\today\quad\militarytime}
          \let\@evenfoot\@oddfoot \overfullrule 3pt
          \let\label=\draftlabel
          \let\marginnote=\draftmarginnote
     \def\@eqnnum{(\theequation)\rlap{\kern\marginparsep\tt\@eqnlabel}%
\global\let\@eqnlabel\@vacuum}  }

\def\numberbysection{\@addtoreset{equation}{section}
          \def\theequation{\thesection.\arabic{equation}}}

\def\underline#1{\relax\ifmmode\@@underline#1\else
          $\@@underline{\hbox{#1}}$\relax\fi}

\def\titlepage{\@restonecolfalse\if@twocolumn\@restonecoltrue\onecolumn
       \else \newpage \fi \thispagestyle{empty}\c@page\z@
          \def\thefootnote{\fnsymbol{footnote}} }

\def\endtitlepage{\if@restonecol\twocolumn \else  \fi
          \def\thefootnote{\arabic{footnote}}
          \setcounter{footnote}{0}}  %\c@footnote\z@ }
\relax

\hybrid

\def\beq{\begin{equation}}
\def\eeq{\end{equation}}
\def\p{\partial}

\def\tauh{\tau_{\hbar}}

\begin{document}

\begin{titlepage}

\title{Dispersionless limit of Hirota equations in some
problems of complex analysis}

\author
{A.Zabrodin\thanks{Institute of Biochemical Physics,
Kosygina str. 4, 117334,
Moscow, Russia and ITEP, 117259, Moscow, Russia}}

\date{April 2001}
\maketitle
\vspace{-9cm}

\centerline{\hfill ITEP/TH-15/01}

\vspace{9cm}
\begin{abstract}

The integrable structure, recently revealed
in some classical problems of the
theory of functions in one complex variable,
is discussed. Given a simply connected
domain in the complex plane,
bounded by a simple analytic curve,
we consider the conformal mapping
problem, the Dirichlet boundary problem, and
to the 2D inverse potential problem associated
with the domain. A remarkable family of real-valued
functionals on the space of such domains is constructed.
Regarded as a function of infinitely many variables,
which are properly defined moments of the domain,
any functional from the family gives
a formal solution to the problems listed above.
These functions are shown to obey an infinite set of
dispersionless Hirota equations. This means that they
are $\tau$-functions of an integrable hierarchy.
The hierarchy is identified with the
dispersionless limit of the 2D Toda chain.
In addition to our previous
studies, we show that, with a more general
definition of the moments, this connection is not
specific for any particular solution to the Hirota equations
but reflects the structure of the hierarchy itself.

\end{abstract}

%\vfill

\end{titlepage}

\section{Introduction}

Recently, it was found \cite{MWZ,WZ,KKMWZ}
that some problems of complex
analysis, such as conformal mapping problem,
Dirichlet boundary problem,
and 2D inverse potential problem, have a hidden integrable
structure. It turns out that there is an infinite hierarchy
of non-linear integrable equations such that variations of
sought-for solutions to the problems listed above under deformation
of given data can be found by solving the hierarchy.
The latter is a certain
scaling limit, usually referred to as dispersionless limit,
of the 2D Toda chain hierarchy.

The conformal mapping problem is to find a conformal map
from a given simply connected domain in the complex plane
to the unit disk. By the Riemann mapping theorem, such a map
does exist and, after imposing some simple
normalization conditions, is unique (see, e.g.,
\cite{GC}). To be specific,
let $\gamma$ be a closed analytic curve in the complex plane
and $D_+$, $D_-$ be respectively
the interior and exterior domains with respect to the curve.
By $w(z)$ we denote the conformal map of $D_-$ to the exterior
of the unit circle normalized in such a way that $\infty$ is
taken to $\infty$ and the derivative at $\infty$ is real positive.

The Dirichlet boundary problem in $D_{-}$ sounds as follows:
given a function $\psi (z)$ on the curve $\gamma$, to find
a harmonic function $f(z,\bar z)$ in $D_{-}$ such that
$f(z, \bar z)=\psi (z)$ for $z\in \gamma$. The solution
is given by the formula
\beq\label{green1}
f(z,\bar z) =-\, \frac{1}{\pi i} \oint_{\gamma}\psi (\zeta )
\p_{\zeta}G(z,\zeta ) d\zeta
\eeq
where $G(z,\zeta )$ is the Green function of the Dirichlet
problem in $D_{-}$. The Green function is
uniquely characterized by the properties \cite{Hille}:
(a) It is regular at infinity and
harmonic everywhere in $D_{-}$ but at the point
$z=\zeta$ where it has the logarithmic singularity:
$G(z, \zeta )=\log |z-\zeta | +O(1)$ as $z \to \zeta$;
(b) It is equal to zero on the boundary
$\p D_{-}=\gamma$, i.e., $G(z,\zeta )=0$
for all $\zeta$ if $z\in \gamma$.

To specify what we mean by the inverse potential problem,
we suppose, in the same setting as above, that the domain
$D_{+}$ is filled with a continuously
distributed background electric charge with uniform density.
As is proved in the potential theory, the potential,
$\Phi (z,\bar z)$, generated by the charge, and its
first derivatives (components of the electric field)
are continuous functions at the boundary, i.e.,
$\Phi^+ =\Phi^-$, $\p_z \Phi^+ =\p_z \Phi^-$, $z\in \gamma$.
Here $\Phi^{\pm}$ is the function $\Phi$ restricted to
$D_{\pm}$. $\Phi^-$ is harmonic everywhere
in $D_{-}$ but at infinity, where it has a logarithmic
singularity, while $\Phi^+$ is
harmonic up to a term proportional to $|z|^2$.
Therefore, one can represent $\Phi^{\pm}$ by Taylor series
(multipole expansions) in a neighbourhood of $0$ and $\infty$
respectively. We suppose,
without loss of generality, that 0 is in $D_{+}$.
The coefficients of the multipole
expansions are harmonic moments of $D_{-}$ and $D_{+}$,
respectively. The inverse potential problem is to restore the
shape of the domain given one of the series, either $\Phi^{+}$ or
$\Phi^{-}$. For definitness, we assume that the
coefficients of $\Phi^{+}$ are given and address the
problem how to find the coefficients of $\Phi^{-}$.
(If one knows {\it both} sets of coefficients,
the curve can be easily reconstructed.)
For a review of inverse problems of potential theory
see, e.g., \cite{invpot} and references therein.

The three problems are closely related to each other.
Given, say, the conformal map $w(z)$, the Green function
can be found explicitly via the formula
\beq\label{green2}
G(z, \zeta )=\log \left | \frac{w(z)-
w(\zeta )}{w(z)\overline{w(\zeta )}-1}\right |
\eeq
where the bar means complex conjugation.
In turn, the holomorphic in $z$
part of the Green function is the
conformal map, though normalized in a different way. To
obtain the map $w(z)$, one should take the holomorphic
part of $G(z, \infty )$.
Concerning the inverse potential problem, we note that
it is the Green function that provides a relation between
$\Phi^+$ and $\Phi^-$. Indeed, let us
modify the potential function and introduce $\tilde \Phi$
which is $\Phi$ with subtracted logarithmic singularity
at infinity, so that $\tilde \Phi^-$ is harmonic in $D_-$
including $\infty$.
Then, according to (\ref{green1}), we can write
\beq\label{green3}
\tilde \Phi^- (z,\bar z)
=-\, \frac{1}{\pi i} \oint_{\gamma}\tilde \Phi^+ (\zeta , \bar \zeta )
\p_{\zeta}G(z,\zeta ) d\zeta
\eeq

Let us summarize the main results of \cite{MWZ,WZ,KKMWZ}.
Let $t_k$, $k\geq 1$,
be the harmonic moments of $D_{-}$, and $t_0$ be the area
of $D_{+}$ (divided by $\pi$):
$t_k =-\,\frac{1}{\pi k}\int_{D_-}z^{-k}d^2z$,
$t_0 =\frac{1}{\pi}\int_{D_{+}}d^2z$.
Note that $t_k$ are in general complex numbers while
$t_0$ is real.
In the subsequent
formulas, the complex conjugate moments, $\bar t_k$,
are considered as independent variables.
Introduce the operators
\beq\label{D1}
D(z)=\sum_{k\geq
1}\frac{z^{-k}}{k}\p_{t_k}\,,
\;\;\;\;\;\;
  \bar D(\bar z)=\sum_{k\geq
1}\frac{\bar z^{-k}}{k}\p_{\bar t_k}
\eeq and
\beq\label{D2}
{\cal D}(z, \bar z)=\p_{t_0}+D(z)+\bar D(\bar z)
\eeq
There exists a real-valued function $F$ of $t_0$, $\{t_k\}$,
$\{\bar t_k\}$ such that the conformal map $w(z)$ is given by
\beq\label{wz}
w(z)=z\exp \left (-\frac{1}{2}\p_{t_0}^{2}F -\p_{t_0}D(z)F\right )
\eeq
This function admits an explicit integral representation.
Moreover, the Green function of the Dirichlet
boundary problem is represented as
\beq\label{green4}
G(z_1 , z_2 )=\log |z_{1}^{-1} -z_{2}^{-1}|+
\frac{1}{2}{\cal D}(z_1 , \bar z_1 )
{\cal D}(z_2 , \bar z_2 )F
\eeq
Let us note that a formula of this type was first
conjectured by L.Takhtajan. See \cite{Takh} for
a rigorous proof and discussion.

The function $F$ is the dispersionless limit of the
logarithm of the $\tau$-function of the 2D Toda chain
hierarchy. It obeys the dispersionless Hirota equations
\beq\label{Hir1} (z_1-z_2 )e^{D(z_1 )D(z_2 )F}
=z_1 e^{-\p_{t_0}D(z_1 )F}
-z_2 e^{-\p_{t_0}D(z_2 )F}
\eeq
\beq\label{Hir2}
z_1 \bar z_2
\left (1-e^{-D(z_1)\bar D(\bar z_2 )F}\right ) =
e^{\p_{t_0}(\p_{t_0}+ D(z_1 )+\bar D(\bar z_2 ))F}
\eeq

The dispersionless Toda hierarchy is an example
of the universal Whitham hierarchy
introduced in \cite{K1}.
It is a multi-dimensional extension of the
hierarchies of hydrodynamic type \cite{hydro,Tsarev}.
Equations of the Toda
hierarchy are known to have infinitely many solutions.
The solutions can be
parametrized \cite{TakTak} by canonical transformations
in a two-dimensional
phase space in such a way that any solution
corresponds to a canonical pair
of functions. (In \cite{TakTak}, the canonical
pair is called twistor data of the solution.)
In fact, this is equivalent to the characterization
of the solutions via string equations.

In the approach developed in \cite{WZ}, it was
one particular solution
to the dispersionless 2D Toda hierarchy that described
conformal maps. The place, if any, of other solutions
in the conformal map context was not clear.
In this paper we address exactly this question.
Extending the previous result, we show that the same class
of conformal maps can be described by {\it any} generic
solution to the hierarchy.

This description is achieved by another,
more general definition
of moments of the domain. The moments are again the
independent flows of the hierarchy. However, their
relation to the shape of the domain is different,
and some additional data is required for their definition.
The additional data is a
real analytic function on the complex plane, $\sigma (z, \bar z)$.
In the context of the inverse potential problem,
this function is density of the background charge,
which was constant in the standard formulation of the
problem. More precisely, we suppose that $\sigma (z, \bar z)$
is given in the whole complex plane and fixed once and for all
while the shape of the domain may vary. The background charge
has density $\sigma (z, \bar z)$ inside the domain and 0 outside.
The inverse potential problem is posed in the same
way as before. The moments, however, are now defined with
respect to the density $\sigma$. In the sequel we refer to
this problem as generalized inverse potential problem.

The result shows that it is the hierarchy of Hirota equations
(rather than particular string equations or other solution-specific
ingredients) which reflects the integrable structure of conformal
maps and the Dirichlet problem.

In Sec.\,2 we introduce the function $F$ with the help
of a variational principle and find its first derivatives.
The first derivatives give a formal solution to the
generalized inverse potential problem.
Formal solutions to the conformal map problem and
the Dirichlet problem are given by second derivatives
of $F$. There is an infinite set of relations
for the second derivatives, which can be combined
into the (dispersionless limit of) Hirota equation.
This stuff is discussed in Sec.\,3.
At last, Sec.\,4 is devoted to an alternative approach
to the dispersionless integrable hierarchy, which allows
one to prove the Lax-Sato equations for conformal maps.

\section{Generalized inverse potential problem}

The inverse potential problem is perhaps
the most convenient context to introduce
the $\tau$-function for curves.
Following the idea of \cite{KKMWZ},
we define the $\tau$-function with the help of an
electrostatic variational principle.
The set-up and the notation are the same as in Sec.\,1.

We recall that the domain
$D_+$ is filled with a background electric charge.
Let us consider non-uniform charge distributions such that
density of the charge
is characterized by a real analytic function $\sigma (z, \bar z)$.
We assume that this function is defined all over the
complex plane, the charge density being
equal to $-\sigma$
inside $D_{+}$ and 0 outside.
With a fixed function $\sigma$,
the 2D electrostatic potential $\Phi$ generated
by the domain is a functional of the shape of the domain.
The potential $\Phi$ obeys the Poisson equation
\beq
\label{1}    - \p_{z}\p_{\bar z}\Phi(z,\bar z)   =
\cases{{\sigma (z, \bar z) } &
   \quad if  \quad
   $z=x+iy \in
D_+$  \cr
\,\,\, 0&  \quad if  \quad  $z=x+iy \in
D_-$\   \cr}
\eeq
The  potential $\Phi$ can be  written as an integral over the
domain $D^+$:
\beq
\label{1bis}
\Phi (z,\bar z)=-\, \frac{2}{\pi}\int_{D_{+}}d^2 z'\
\sigma (z', \bar z') \log|z-z'|
\eeq

In $D_{-}$, the potential is a harmonic function
represented by the series
\beq\label{2}\Phi^{-}(z,\bar z) =-2t_{0}\log|z|+
2{\cal R}e\sum_{k>0}\frac{v_k}{k} z^{-k},
\eeq
as $z\to \infty$, where
\beq\label{3}
v_k=\frac{1}{\pi}\int_{D_+} z^k \sigma (z,\bar z)d^2 z
    \qquad (k > 0)
\eeq
are moments of
the interior domain $D_{+}$
with respect to the density $\sigma$, and
\beq\label{3bis}
\pi t_0 =\int_{D_{+}}\sigma (z, \bar z)d^2z
\eeq
is the total charge of $D_{+}$.

In $D_+$, the potential (\ref{1bis}) is represented
by a series $\Phi^+$ which
has both harmonic and non-harmonic terms.
The latter is due to the background charge.
We denote the non-harmonic part by $U(z,\bar z)$.
The function $U$ obeys the Poisson equation
$\p_z \p_{\bar z} U(z,\bar z)=\sigma (z,\bar z)$ and is
fixed by the condition
that its expansion around $0$ (recall that
$0\in D_{+}$) does not contain harmonic terms.
For example, if $\sigma (z, \bar z)=1$ then $U(z,\bar z)=z\bar z$.
We adopt the following general form of $U$:
\beq\label{U}
U(z,\bar z )=- \sum_{m,n\geq 1}T_{mn}z^{m}\bar z^n
\eeq

The expansion of $\Phi$ around $z=0$
can be obtained from (\ref{1bis}) in a direct way.
The result is:
\beq\label{4}
\Phi^{+}(z,\bar z)= - U(z, \bar z)+ v_0+ 2{\cal R}e
\sum_{k>0}t_k z^k
\eeq
Here the non-harmonic term, $U$, is determined by
the background charge, as is explained above,
and does not depend on the form of the domain.
It is the harmonic contribution that takes the boundary
into account. The coefficients are
\beq\label{5}
t_k= \frac{1}{2\pi i k}\oint_{\gamma}z^{-k}
\p_z U \,d z \qquad (k > 0)
\eeq
They provide a complimentary (to $v_k$)
set of parameters characterizing the domain $D_{+}$.
Here and below we use the same notation for the moments
($t_0$, $t_k$, etc) as in Sec.\,1.
In the case $\sigma =1$ they are equal to harmonic
moments of the exterior domain $D_{-}$. Finally,
\beq
\label{v0}
v_0=\frac{2}{\pi}\int_{D_+}  \log|z|\, \sigma (z, \bar z)\,d^2z.
\eeq
Note that the moments $v_k$ (\ref{3}) and $t_0$ (\ref{3bis})
also admit a similar contour integral representation:
\beq\label{3a}
v_k=\frac{1}{2\pi i}\oint_{\gamma} z^k \p_z U \,dz
    \qquad (k > 0)\,,
\;\;\;\;\;
t_0=\frac{1}{2\pi i}\oint_{\gamma}\p_z U \,dz
\eeq
which is a simple consequence of the Green formula.

Since the potential and the electric field
are continuous at the boundary of the
domain, the two complementary sets of moments
are related by the conditions
\beq\label{cont}
\Phi^+=\Phi^-\,, \;\;\;\;
\p_{z} \Phi^+ =\p_{z} \Phi^-\,,
\;\;\;\;\;
z\in \gamma
\eeq
The generalized inverse potential problem is to find
the curve $\gamma$ given the function $\sigma$ and
one of the functions $\Phi^+$ or $\Phi^-$,
i.e., one of the infinite sets of moments.
We will choose as independent
variables the total charge $\pi t_0$ and the
coefficients $t_k , \,  \bar t_k \ (k\ge 1)$ of the
multipole expansion of $\Phi$ around the origin. Under certain
conditions, which we do not discuss here, they uniquely determine
the form of the curve, at least if the deformed curve to be
described by $t_k$ is close enough to a fixed given curve.
(Some important results concerning uniqueness of solution to
inverse potential problems can be found in \cite{uniq}.)
One such condition, specific for the case of non-uniform
density of the background charge, is worth mentioning.
For the problem to be well-posed,
one should require that $\sigma (z, \bar z)\neq 0$
at any point of the curve $\gamma$. Otherwise
a small perturbation of the curve around the zero
of $\sigma$ does not result, at least in the first
order, in the variation of $\Phi$. To avoid this difficulty,
one may assume, for example, that $\sigma (z, \bar z)>0$.

In other words, under the conditions mentioned above,
$\{t_k\}_{k=0}^{\infty}$ is a good set
of local coordinates
in the space of analytic curves.
The moments $v_j$ are also determined by
$\{t_k\}_{k=0}^{\infty}$, so any moment of the interior
can be regarded as a function of $t_k$.
For simplicity,
we assume in this paper that
only a finite number of $t_k$ and $T_{mn}$ are
non-zero. In this case the series (\ref{4}) is a polynomial in $z$,
$\bar z$ and, therefore, it defines the function $\Phi^{+}$ for
all $z\in D_{+}$.
As soon as $\sigma$ is real, $t_0$, $v_0$ are real
numbers while all other
moments are in general complex quantities.

We now pass to the electrostatic variational principle
mentioned in the beginning of this section.
Consider the energy functional
describing a
charge with a density $\rho (z,\bar z)$ in the background
potential $\Phi$ (\ref{1bis})
generated by the charge distributed in $D_{+}$ with the
density $\sigma$:
\beq\label{32B}
{\cal
E}\{\rho\}=-\frac{1}{\pi^2}
\int\! \int d^2 z \,d^2 z' \  \rho(z,
\bar z)\  \log|z-z'|\rho(z',\bar
z')\ \,-\frac{1}{\pi}\int\! d^2
z\ \rho(z,\bar z) \ \Phi(z,\bar z).
\eeq
The first term is the 2D
``Coulomb'' energy of the charge
while the second one is the energy
of interaction with the background charge.
Clearly, the distribution of the
charge neutralizing the background charge
gives minimum to the functional.
We denote density of this distribution by $\rho_0$:
$\rho_{0}=\sigma$ inside the domain and $\rho_{0}=0$
outside, so that the total charge density be zero.
At the minimum the functional is equal to minus
electrostatic energy
$E$ of the
background charge $\mbox{min}_{\rho}{\cal E}(\rho )=-E$, where
\beq
\label{32A}
E= -\, \frac{1}{\pi^2}
\int \!\!\int_{D_{+}} \!\! d^2 z  d^2 z' \
\sigma (z, \bar z)\log|z-z'|
\sigma (z', \bar z')
\eeq
Varying $\rho$ in (\ref{32B}) and then setting
$\rho=\rho_0$, we
obtain eq.\,(\ref{4}).

Let us find derivatives of $E$
with respect to $t_k$.
This can be done in different ways.
One of them would be to take the derivative
of (\ref{32A}) directly. In doing so, one should
somehow control the contribution coming from
an infinitesimal change of the domain under the
variation of $t_k$. To avoid this problem, we
use the same argument as in the case
$\sigma =1$ (see \cite{KKMWZ}). Namely,
instead of taking the derivative of (\ref{32A}),
we differentiate
$E$ represented as $-{\cal E}(\rho )$
in the form (\ref{32B}) at the extremum.
There are contributions of two kinds: one from the explicit
$t_k$-dependence of $\Phi^+$ and another from
the implicit $t_k$-dependence of $\rho_0$:
$$
\left. \frac{\p E}{\p t_k} =
-\,\frac{\p {\cal E}(\rho )}{\p t_k}\right |_{\rho =\rho_0}=
\left.
\frac{1}{\pi}\int_{D_{+}}\!\! d^2 z \sigma (z,\bar z)\p_{t_k}
\! \Phi^+ (z,\bar z)
-\,\int \! d^2 z\,
\frac{\delta {\cal E}(\rho )}{\delta
\rho (z,\bar z)}\right |_{\rho =\rho_0}
\!\!\!\!
\frac{\delta \rho_0 (z,\bar z)}{\delta t_k}
$$
Since we are at the extremum, the
variational derivative in the second term is zero,
and this term does not contribute. In other words,
the derivative is equal to the partial derivative
of $-{\cal E}$ calculated at the extremum.
Let us treat, for a while,
$v_0$ and $t_k$ as independent variables, then
the partial derivatives of $\Phi^+$ are especially
simple:
$\p_{t_k}\Phi^+ (z,\bar z)=z^k$,
$\p_{v_0}\Phi^+ (z,\bar z)=1$.
Plugging them into the integral,
we arrive at the relations
\beq\label{25a}
\frac{\p E}{\p t_k}= v_k, \;\;\;
\frac{\p E}{\p
\bar t_k}= \bar v_k,
\;\;\;
\frac{\p E}{\p v_0}=- t_0,
\eeq
where the
partial derivative with respect to
$t_k$ is taken at fixed $v_0$ and
$t_j \ (j   \neq 0, k)$.
Therefore, the differential $dE$ is
$dE= \sum_{k>0}(v_k dt_k +\bar v_k d\bar t_k )
- t_0 dv_0$. Note that there are two
contributions to the $t_k$-derivative
of (\ref{32A}): one from $\p_{t_k}\Phi$ and another
from the change of the domain. The former is easily
calculated to be $\frac{1}{2}v_k$, thus the
latter is $\frac{1}{2}v_k$, too. However, it is
not so easy to obtain this result directly (see the
appendix for an example of such a calculation).
A direct proof of the formulas for first derivatives
in the case $\sigma =1$ can be found in \cite{Takh}.

It is more natural, however, to treat the total charge
rather than $v_0$ as
an independent variable, i.e., to apply the
variational principle
at a fixed total charge. This
is achieved via the
Legendre transformation. Let us introduce the function
$F=E+ t_0 v_0$, whose  differential is
$dF= \sum_{k>0}(v_k dt_k
+\bar v_k d\bar t_k ) + v_0 dt_0$.
The integral representation of this function follows from
(\ref{32A}):
\beq\label{31}
F=-\, \frac{1}{\pi^2}
\int \!\!\int_{D_{+}} \!\! d^2 z  d^2 z' \
\sigma (z, \bar z)\log|z^{-1}-z'^{-1}|
\sigma (z', \bar z')
\eeq

The function $F$ plays a major role
in what follows. It is a real-valued function
of the moments $t_0 , \, t_1 , \, t_2 , \ldots $.
Rewriting (\ref{25a}) in the new variables,
we get the main property of the
function $F$:
\beq\label{25}
\frac{\p F}{\p t_k}=v_k, \;\;\:
\frac{\p F}{\p \bar t_k}=\bar v_k,
\;\;\;
\frac{\p F}{\p t_0}=v_0
\eeq
where the derivative with respect to $t_k$ is taken at
fixed $t_j$ ($j \neq k$).
The very existence of the common potential function
for moments implies
the symmetry relations for their derivatives,
$\p_{t_n} v_k=\p_{t_k} v_n$,
$\p_{\bar t_n} v_k=\p_{t_k} \bar v_n$,
which were first obtained in \cite{MWZ} (for the case
$\sigma =1$) within a different approach.

To represent the result (\ref{25}) in a more compact form,
let us modify the potential $\Phi$ and introduce
$\tilde\Phi(z,\bar z)=\Phi (z, \bar z)
+2 t_0 \, \mbox{log}\,|z|+v_0$, so that
\beq\label{tau1}
\tilde\Phi(z,\bar z)=
-\, \frac{2}{\pi}\int_{D_{+}}\! \!d^2 \zeta\
\sigma (\zeta ,\bar \zeta )
\log \left |z^{-1}- \zeta^{-1}\right |
\eeq
and
\beq\label{31a}
F= \frac{1}{2\pi}
\int_{D_{+}} \!\!\!\! d^2 z \ \sigma (z, \bar z)\,\tilde
\Phi^{+}(z,\bar z)
\eeq
The modification amounts to adding the neutralizing
point-like charge at the origin.
The expansion of the modified potential around the origin is
\beq\label{or}
\tilde\Phi^+ (z, \bar z)=-U(z,\bar z)+
2t_0\log |z|+\sum_{k\geq 1}(t_k z^k +\bar t_k \bar z^k)
\eeq
From now on, this formula will serve as a definition
of the variables $t_k$.
Likewise $\Phi$, the function $\tilde \Phi$ and its first
derivatives are continuous on
the boundary $\gamma =\p D_{+}$.
For $z$ in $D_{-}$ the function $\tilde \Phi (z,\bar z)$
is harmonic.
In terms of $\tilde\Phi$,
formulas (\ref{25}) can be compactly written in the form
\beq\label{fir2}
{\cal D}(z,\bar z)F=\tilde \Phi^{-} (z,\bar z)
\eeq
where we use the operator ${\cal D}$ introduced in (\ref{D2}),
and $z$ is assumed to be in $D_{-}$.
Both sides are to be understood as series in $z$ and $\bar z$.
It is our assumption that they converge in $D_{-}$.

To conclude this section, let us discuss homogeneity
properties of the function $F$.
Under the assumption that only a finite number
of $t_k$'s are non-zero,
we can substitute (\ref{4}) into (\ref{31})
and perform the term-wise integration.
We get the following relation:
\beq
\label{free2}
2F = -\frac{1}{4\pi}\int_{D_{+}}\!\!d^2 z U\Delta U
+t_0 v_0 +
\sum_{k>0}(t_k v_k +\bar t_k \bar v_k)\,,
\eeq
where $\Delta =4\p_z \p_{\bar z}$ is the Laplace operator.
Generally speaking, (\ref{free2}) does not have
any definite homogeneity properties as a
function of $t_k$. However, for homogeneous functions
$\sigma$, $F$ is quasihomogeneous. For instance,
if $\sigma (z, \bar z)=|z|^{2M-2}$ with an integer $M$,
then (\ref{free2}) can be brought into the form
$$
4MF= -t_{0}^{2} +2M
t_0 \p_{t_0}F
+\sum_{k\geq 1}(2M-k)(
t_k \p_{t_k}F + \bar t_k \p_{\bar t_k}F )
$$
where the relations (\ref{25}) are taken into account.
This relation reflects homoheneity
properties of the moments under the
rescaling $z \to \lambda z$ and
means that $F$ is a quasihomogeneous
function of degree $4M$.

A more fundamental homogeneity property,
valid for any $\sigma$,
can be revealed by extending the definition
of the function $F$. Let us treat parameters
$T_{mn}$ defining $U$ or $\sigma$ (see (\ref{U}))
on equal footing with the variables $t_k$.
To wit, let us identify
$t_k = T_{k0}$,
$\bar t_k = T_{0k}$,
$t_0 = T_{00}$, and allow $F$ to be a function
of all $T_{mn}$, $m,n \geq 0$, regarded as an
extended set of independent variables.
Using the same variational argument as before,
one can deduce the relations
$$
\frac{\p F}{\p T_{mn}}=\frac{1}{\pi}
\int_{D_{+}} z^n \bar z^m \sigma (z, \bar z)d^2z
$$
which generalize (\ref{25}). Then (\ref{free2})
acquires the form of the pure homogeneity condition:
\beq\label{hom}
2F = \sum_{n,m\geq 0} T_{mn}\frac{\p F}{\p T_{mn}}
\eeq

The $\tau$-function
can be defined as $ \tau =e^F$.
This notation, however, is introduced here
only for the purpose to stress the relation
to Hirota equations.
It is $F$, i.e., the logarithm of the
$\tau$-function, which we only deal with in the sequel.
Recall that the $\tau$-function itself does not have a good
dispersionless limit. Only its logarithm, multiplied
by a small dispersion parameter, makes
sense in this limit. The bilinear
Hirota equations for $\tau$ are rewritten
as highly non-linear equations for $\log\tau$.
In the next section, we derive these equations
starting from the definition of the function $F$.

\section{Dirichlet problem
and dispersionless Hirota equations for the
function $F$}

The dispersionless Hirota equations are certain
relations between second derivatives
of logarithm of the $\tau$-function with respect
to $t_k$. The second derivatives of $F$
are coefficients of the expressions like $D(z_1)D(z_2)F$,
where the operator $D(z)$ is introduced in (\ref{D1}).
Recall that
we already know first derivatives of $F$,
see (\ref{fir2}). Therefore, what we need to find
is ${\cal D}(z_1, \bar z_1)\tilde \Phi (z_2 , \bar z_2)$
for $z_1 , z_2 \in D_{-}$.
For this we use the
following general argument.

Consider a small
deformation of the domain $D_{+}$ obtained from the original one
by adding to it an infinitesimal bump (of
arbitrary form) with area $\epsilon$
located at a point $\xi \in \gamma$.
By $\delta_{\epsilon (\xi )}$, denote variation
of any quantity under such deformation in the first
order in $\epsilon$.
The potential generated by the deformed domain
is then $\tilde\Phi +\delta_{\epsilon (\xi )}\tilde\Phi$,
where
\beq\label{or1}
\delta_{\epsilon (\xi )}\tilde\Phi (z, \bar z)=
-\,\frac{\epsilon}{\pi}\sigma (\xi , \bar \xi )
\log |z^{-1}-\xi^{-1}|^2
\eeq
is the potential generated by the bump, as it is clear
from the explicit formula (\ref{tau1}).
Comparing the expansion of
$\delta_{\epsilon (\xi )}\tilde\Phi (z, \bar z)$
in $z$ as $z\to 0$ with (\ref{or}), we find
variations of $t_k$ under the deformation:
$$
\delta_{\epsilon (\xi )}t_0 =\frac{\epsilon}{\pi}
\sigma (\xi ,\bar \xi )\,,
\;\;\;\;\;
\delta_{\epsilon (\xi )}t_k =\frac{\epsilon}{\pi k}
\sigma (\xi ,\bar \xi )\xi^{-k}\,,
\;\;\;\;\;k\geq 1
$$
Given any function $A$ of the moments, its variation,
$\delta_{\epsilon (\xi)}A$, is
$$\delta_{\epsilon (\xi)}A = \p_{t_0}\! A
\delta_{\epsilon (\xi)}t_0
+\sum_{k\geq 1}(
\p_{t_k}\! A
\delta_{\epsilon (\xi)}t_k +
\p_{\bar t_k}\! A
\delta_{\epsilon (\xi)}\bar t_k )
$$
which can be written as
\beq\label{D4}\delta_{\epsilon(\xi)} A =
\frac{\epsilon}{\pi} \sigma (\xi , \bar \xi )
{\cal D}(\xi, \bar \xi )A\,,
\;\;\;\;\; \xi \in \gamma
\eeq
So we see that for
$\xi$ on the boundary curve $\gamma =\p D_{+}$
the operator ${\cal D}(\xi, \bar \xi)$ has a clear
geometrical meaning: the result of its action on any
quantity is proportional to the variation of this
quantity under attaching the bump at the point $\xi$.
To put it differently, we can say that the
boundary value of the function
${\cal D}(z,\bar z)A$ is equal to
$\pi \delta_{\epsilon (z)}\! A/(\epsilon \sigma (z,\bar z))$,
$z\in \gamma$. For functions $A$ such that the series
$D(z)A$ converges everywhere in $D_{-}$ up to the
boundary,
this remark gives a usable method to find the function
${\cal D}(z,\bar z)A$. To wit, this function is
harmonic in $D_{-}$, and its boundary value is given
by (\ref{D4}). It is then just the subject of the Dirichlet
boundary problem to find the harmonic function given
its boundary value. This method appears to be especially
useful when one is able to find the left hand side of (\ref{D4})
independently.

As a very simple example, one can easily
find $\delta_{\epsilon (\xi)}F =
\frac{\epsilon}{\pi}\sigma (\xi , \bar \xi )
\tilde \Phi (\xi , \bar \xi )$ from the integral
representation of $\tilde\Phi$. In this case,
the harmonic continuation does not require
any additional care since the function
$\tilde\Phi$ is already harmonic in $D_{-}$.
Whence we reproduce the result (\ref{fir2}).

Let us turn to second derivatives. As it follows from
(\ref{D4}), the boundary
value of the function
${\cal D}(\xi , \bar \xi ) {\cal D}(z , \bar z )F=
{\cal D}(\xi , \bar \xi ) \tilde\Phi (z , \bar z )$,
regarded as a function of $\xi$, is given by eq.\,(\ref{or1})
for $z$ in $D_{-}$:
\beq\label{sec1}
{\cal D}(\xi , \bar \xi ) {\cal D}(z , \bar z )F=
-2\log |z^{-1}-\xi^{-1}|\,,
\;\;\;\;\;\; \xi \in \gamma
\eeq
Here we do not discuss convergence of the
series $D(\xi )\tilde\Phi$ in $D_{-}$, which is
the necessary assumption
for the harmonic continuation rule explained above.
To be sure that the class of domains having this property
is not empty, one may think of domains
and background charges that are close enough
to the uniformly charged disk,
in which case the proof of the convergence
is the matter of some routine estimates.
The right hand side of (\ref{sec1})
is not harmonic in $D_{-}$ as it stands
because of the logarithmic singularity at $\xi =z$.
To find the harmonic continuation of this function,
we add to (\ref{sec1}) a function harmonic everywhere
in $D_{-}$ but at the point $z$, where it has the logarithmic
singularity $+2\log |z^{-1}-\xi^{-1}|$ (thus cancelling
the singularity of (\ref{sec1})), and equal to zero on the
boundary. By definition of the Green function of the
Dirichlet problem (see Sec.\,1),
it is the Green function $G(z, \xi )$ (\ref{green2}).
Hence we obtain formula (\ref{green4}).
Tending $z_2$ to $\infty$ and then separating the
holomorphic part in $z_1$, one arrives at the expression
(\ref{wz}) for the conformal map.
Separating holomorphic and antiholomorphic parts
of (\ref{green4}) in both
variables, and using (\ref{green2}), we obtain:
\beq\label{sec5}
\log\frac{w(z_1)-w(z_2)}{z_1 -z_2} =-\,
\frac{1}{2}\p_{t_0}^{2}F +D(z_1)D(z_2)F
\eeq
\beq\label{511}
\log \left ( 1-\frac{1}{w(z_1)\overline{w ( z_2) }}\right )
=\,- D(z_1)D(z_2)F
\eeq
The limiting case of
(\ref{sec5}) as $z_2 \to \infty$ gives another
useful formula for the conformal map:
\beq\label{wza}
w(z)=e^{-\frac{1}{2}\p_{t_0}^{2}F }\left (
z- (\p_{t_0} +D(z) )\p_{t_1}F \right )
\eeq
All these formulas have the same form for any
function $\sigma$.

It is now straightforward to see how the Hirota
equations come into play. They are obtained by
plugging (\ref{wz}) into equalities of the type
(\ref{sec5}), (\ref{511}), thus excluding $w(z)$
from them. In this way one obtains eqs.\,(\ref{Hir1}),
(\ref{Hir2}), in which one recognizes
dispersionless limit of bilinear Hirota equations
for the 2D Toda chain hierarchy.
Another Hirota equation does not even require
(\ref{wz}) for its derivation: take
three copies of eq.\,(\ref{sec5}) for each pair
of the points $z_1 , z_2 , z_3$, take the
exp-function of both sides of each equation and then
sum them up. One arrives at the equation
\beq\label{KP1}
(z_1 -z_2)e^{D(z_1)D(z_2)F}+
(z_2 -z_3)e^{D(z_2)D(z_3)F}+
(z_3 -z_1)e^{D(z_3)D(z_1)F} =0
\eeq
This is the dispersionless KP hierarchy (which is
a holomorphic, with respect to the variables $t_k$,
sector of the Toda chain)
written in the most symmetric form. Another version, known in
the literature \cite{CK,TakTak},
is obtained from here as $z_3 \to \infty$.

The three equations,
(\ref{Hir1}), (\ref{Hir2}), and
(\ref{KP1}) (along with their bar-versions)
are to be understood as an infinite set of
relations for second derivatives of $F$.
Altogether, they form the integrable hierarchy.
These relations are obtained by expanding
both sides in a power series in $z_i$ and
comparing the coefficients. For example,
comparing leading terms in both sides
of (\ref{Hir2}) as $z_1, z_2 \to \infty$, one
obtaines the first equation of the hierarchy.
It is the dispersionless version of the Toda equation:
$\p_{t_1 \bar t_1}^{2}F=e^{\p_{t_0}^{2}F}$.
(In terms of the field $\phi =v_0 =\p_{t_0}F$,
it is
$\p_{t_1 \bar t_1}^{2}\phi =\p_{t_0} e^{\p_{t_0}\phi }$.)
Next-to-leading terms give the equation
$\p_{t_2 \bar t_1}^{2}F= 2
\p_{t_0 t_1}^{2}F\,
\p_{t_1 \bar t_1}^{2}F$.

An important comment is in order.
All the formulas for second derivatives
of $F$, obtained in this section, are of the same
form for any function $\sigma$, i.e., $\sigma$ does not
enter them explicitly. Indeed, consider, for instance,
(\ref{sec5}), (\ref{511}). Their left hand sides, being
expressed through the conformal map only,
are determined mearly by the shape of the domain.
For a given domain, they are the same
for any density of the background charge.
Clearly, this also holds true for the generating
formula for second derivatives (\ref{green4}).
Note that $F$ itself does
depend on the choice of $\sigma$.

This remark can be reformulated as a covariance property
of the second derivatives.
Let $\hat F$, $\hat t_k$ be the dependent and independent
variables related to another density function,
$\hat \sigma (z, \bar z)$.
As it follows from the above comment,
second derivatives of $F$ are covariant
in the following sense:
\beq\label{hat}
\frac{\p^2 \hat F}{\p \hat t_j \p \hat t_k}=
\frac{\p^2 F}{\p t_j \p t_k}
\eeq
when the both sides are calculated for the
same domain $D_{+}$.
A similar formula for mixed
$t_i , \bar t_k$-derivatives also holds true.
We stress that this covariance takes place only for
second derivatives. For instance, as it is seen
from (\ref{fir2}), in general
$\p_{\hat t_j}\hat F \neq \p_{t_j}F$, since the potential
$\tilde\Phi^{-}$ does depend on $\sigma$ for a given domain.

The Hirota equations are invariant with respect to
varying the density $\sigma$. At the same time,
the particular solution describing the conformal maps
and the Dirichlet problem is determined by the background
charge. We conclude that the function $\sigma$
parametrizes different solutions of the Hirota equations.
We come back to this point again, within
a different approach, in Sec.\,5.
Regarded as a function of all $T_{mn}$ (see the
end of Sec.\,2), $F$ solves the extended Toda chain
hierarchy introduced in \cite{BX} in the context of
the 2-matrix model.

A short comment on the bilinear
Hirota equations and their dispersionless limit is in order.
Let $\hbar$ be a parameter.
The equations are written for a
$\tau$-function $\tauh$ which depends on the times
$t_k$ and on the parameter $\hbar$.
The list of equations is parallel to the one
for the dispersionless case. We use the same notation.
Eqs.\,(\ref{Hir1}), (\ref{Hir2})
are replaced by
\begin{eqnarray}\label{fHir3}
&&z_1\left (e^{\hbar (\p_{t_0} -D(z_1))}\tauh \right )
\left (e^{-\hbar D(z_2 )}\tauh \right )
- z_2 \left (e^{\hbar (\p_{t_0} -D(z_2 ))}\tauh \right )
\left (e^{-\hbar D(z_1 )}\tauh \right )\,=
\nonumber\\
&=&(z_1-z_2 )\left (e^{-\hbar (D(z_1)+D(z_2 ))}\tauh \right )
\left (e^{\hbar \p_{t_0} }\tauh \right )
\end{eqnarray}
\begin{eqnarray}\label{fHir4}
&&\left (e^{-\hbar D(z_1)}\tauh \right )\!
\left (e^{-\hbar \bar D(\bar z_2 )}\tauh \right )
-\tauh
\!\left (e^{\hbar (\bar D(\bar z_2 )
-D(z_1 ))}\tauh \right )\! =
\nonumber\\
&=&(z_1\bar z_2 )^{-1}\!
\left (e^{-\hbar (\p_{t_0} +D(z_1))}\tauh \right )
\!\left (e^{\hbar (\p_{t_0} +\bar D(\bar z_2 ))}\tauh \right )
\end{eqnarray}
respectively. The dispersionful analog of eq.\,(\ref{KP1}) is
\beq\label{fHir1}
(z_1 -z_2)\left (e^{-\hbar (D(z_1)+D(z_2 ))}\tauh \right )
\left ( e^{-\hbar D(z_3 )}\tauh \right )+[\mbox{cyclic per-s of}
\,\, z_1, z_2 ,z_3 ]=0
\eeq
The operators $e^{\hbar D(z)}$ are shift operators,
so the above equations are difference ones.
To make this explicit, introduce the vector fields
$\p_{x_i}=-D(z_i )$.
In the variables $x_i$ the equations aquire a more familiar
form. For example, eq.\,(\ref{fHir1}) reads
\beq\label{fHir1a}
(z_1 -z_2)\tauh (x_1 +\hbar ,\, x_2 +\hbar , \, x_3)
\tauh (x_1 ,\, x_2 , \, x_3 +\hbar )
+[\mbox{cyclic per-s of}
\,\, 1, 2 , 3 ]=0
\eeq
These bilinear difference equations
for the $\tau$-function enjoy many remarkable
properties.
Eq.\,(\ref{fHir1}) first appeared in Hirota's paper \cite{Hirota}.
For a review of the difference Hirota equations see e.g.\,\cite{Z}.

As is clear from these formulas, the parameter $\hbar$
plays the role of lattice spacing. In the
dispersionless limit, difference equations become
differential ones. The limit is well-defined
on the class of solutions such that there
exists a finite limit
\beq\label{ful1}
F=\lim_{\hbar \to 0} \hbar^2 \log\tauh
\eeq
Rewriting the Hirota equations for $\log\tauh$ and
extracting the leading terms as $\hbar \to 0$, one arrives
at the dispersionless Hirota equations for $F$.
Different aspects of the mathematical theory of
dispersionless integrable equations were discussed
in \cite{GiKo1,Dubr,TakTak,CK}.

\section{Integrable structure
of conformal maps in the Ta\-ka\-sa\-ki-Ta\-ke\-be formalism}

The approach of the previuous section
is maybe the easiest way to see the link between Dirichlet problem
and Hirota equations.
However, such more customary attributes of integrability
like Lax representation remain obscure, from this point
of view.
In this section we discuss an alternative approach which makes
the integrable structure explicit via the Lax-type
equations for the inverse conformal map, $z(w)$.
In our normalization,
the general form of the inverse map is
\beq
\label{z}
z(w)=r w+\sum_{j=0}^{\infty}u_j w^{-j}.
\eeq
where $r$ is real positive. For $w$ on
the unit circle, eq.(\ref{z}) gives a
parametrization of the curve.

We introduce one more new object.
Let $S(z)$ be the analytic continuation
of the function $\p_z U$ away from the curve:
\beq\label{sch1}
S(z) =\p_z U(z, \bar z),\,
\;\;\;\;\;z\in \gamma
\eeq
For analytic curves, $S(z)$ can be proved
to be holomorphic in some strip-like neighbourhood
of the curve $\gamma$.
This function turns out to be a very useful
object from technical point of view.
We call it {\it generalized Schwarz function},
or simply Schwarz function for
brevity. (The latter name is commonly used in the case
$\sigma =1$, see e.g. the book \cite{Davies}.)
Comparing its definition with the definition
of the moments through contour integrals
(\ref{5}), (\ref{3a}), one obtains the
Laurent series representation of the Schwarz
function:
\beq
\label{sch2}
S(z) =\sum_{k=1}^\infty k
t_{k}z^{k-1}+\frac{t_0}{z}+\sum_{k=1}^\infty v_k z^{-k-1}.
\eeq
We thus see that $S(z)$ is the generating function of
{\it all} moments $t_k , v_k$. We regard it as a function
of the moments $t_k$: $S(z)=S(t_k , z)$.
Given the coefficients $T_{mn}$ (\ref{U})
and the Schwarz function, the curve $\gamma$
is determined by eq.\,(\ref{sch1}).

Following the lines
of ref.\,\cite{TakTak}, we introduce the
dispersionless 2D Toda hierarchy starting with the following
data.
Suppose we are given with four Laurent series of the form
\beq
\label{L}
L=rw+\sum_{j=0}^{\infty}u_jw^{-j}\,,
\;\;\;\;\;
M=\sum_{k=1}^{\infty}kt_kL^k +t_0 +
\sum_{k=1}^{\infty} v_k L^{-k}
\eeq
\beq\label{barL}
\bar L=rw^{-1}+\sum_{j=0}^{\infty}\bar u_jw^{j}\,,
\;\;\;\;\;
\bar M=\sum_{k=1}^{\infty}k\bar t_k\bar L^k +t_0 +
\sum_{k=1}^{\infty} \bar v_k \bar L^{-k}
\eeq
The coefficients $t_0, t_k, \bar t_k$ are regarded as
{\it independent variables} while all other
coefficients ($r, u_j, \bar u_j, v_j, \bar v_j$) are
{\it dependent variables}. Imposing some conditions
on the Laurent series, we are going to study the latter
as functions of the former.

Recall now the Takasaki-Takebe theorem \cite{TakTak}.
Let the Poisson bracket $\{,\} $ be defined as
\beq
\label{PB}
\{f,\,g\}=w\frac{\p f}{\p w}\frac{\p g}{\p t_{0}}
-w\frac{\p g}{\p w}\frac{\p f}{\p t_{0} }
\eeq
and let $f,g$ be a canonical pair,
$\{f(w, t_0),\,g(w, t_0)\}=f(w,t_0)$.
Suppose the four Laurent series
$L, \bar L, M, \bar M$ introduced above
are connected by the functional relations
\beq
\label{RH}
\bar L=f^{-1}(L,M)\,,
\;\;\;\;\;
\bar M=g(L,M)
\eeq
Then the following is true:
\begin{itemize}
\item[(A)] The pairs $L,M$ and $(\bar L)^{-1}, \bar M$ are
canonical: $\{L,M\}=L$,
$\{\bar L^{-1},\bar M\}=\bar L^{-1}$;
\item[(B)] The four Laurent series satisfy
the Lax-Sato equations
\beq
\label{B}
\frac{\p X}{\p t_n}=\{H_n , X\}\,,
\;\;\;\;\;\;\;
\frac{\p X}{\p \bar t_n}=-\{\bar H_n , X\}
\eeq
where $X$ stands for any one of $L, \bar L, M, \bar M$ and
\beq
\label{ham1}
H_n=(L^n)_{\geq 1}+\frac{1}{2}(L^n)_{0}\,,
\;\;\;\;\;
\bar H_n=(\bar L^n)_{\leq -1}+\frac{1}{2}(\bar L^n)_{0}\,,
\;\;\;\;\; n\geq 1
\eeq
\item[(C)] There exists a function $F$ of
$t_k$, $\bar t_k$, and $t_0$ such that
$v_k =\p_{t_k}F$, $\bar v_k =\p_{\bar t_k}F$,
and the function $F$ obeys the dispersionless Hirota
equations.
\end{itemize}

This theorem is proved in \cite{TakTak}.
The proof of (A)
consists in differentiating (\ref{RH}) with respect to
$w,t_0$ and combining the results. Taking derivatives
of (\ref{RH}) with respect to $t_n, \bar t_n$ and
using (A), we obtain (B), where
$(...)_{\geq 1}$ ($(...)_{\leq -1}$)
means the part of the Laurent series
with strictly positive (strictly negative)
powers of $w$, and the
constant term is denoted by $(...)_0$.
To prove (C), one should
take $\p M /\p t_k$ at fixed $w$. Using the Lax-Sato
equations (\ref{B}), one gets $\p_{t_j} v_k =
(wL^k dH_j )_0$.
Further, this representation, together with (\ref{ham1}),
implies the symmetry
$\p_{t_j}v_k = \p_{t_k}v_j$, which, in its turn, implies
the existence of a common potential function for $v_k$.
The proof of the Hirota equations for $F$
from these data is implicit
in \cite{TakTak}. More explicit arguments
are given in \cite{CK} (see also \cite{KKMWZ}).

It then follows that $L$ and $\bar L$ are the two Lax
functions of the dispersionless
Toda hierarchy. Equations of the hierarchy
are written for the ``potentials" $r$, $\{u_n\}$,
$\{\bar u_n\}$ as functions of the times $t_0$,
$\{t_n\}$, $\{\bar t_n\}$. These times appear expicitly
in the Laurent series $M$, $\bar M$ which
are the Orlov-Shulman functions \cite{OS} of the hierarchy.
The Takasaki-Takebe theorem says that solutions to the
hierarchy are in one-to-one correspondence with canonical pairs
of functions $f,g$.

Remarkably, this formalism turned out to be
very convenient for describing the integrable
structure of conformal maps.
The identification of the objects related to
conformal maps with the above Laurent series
goes as follows:
$$
L(w)=z(w)\,,
\;\;\;
\bar L(w)=\bar z(w^{-1})\,,
\;\;\;
M(L)=zS(z)\,,
\;\;\;
\bar M(\bar L)=\bar z \bar S(\bar z)
$$
(for a series $f(z)=\sum f_j z^j$ we write
$\bar f (z)=\sum \bar f_j z^j$).
Together with
$S(z)=\p_z U(z,\bar z)$
this means that we have the following
functional relations between $L, \bar L, M, \bar M$:
\beq \label{canon}
M=L\, \displaystyle{ \frac{\p U(L, \bar L)}{\p L}}\,,
\;\;\;\;\;\;
\bar M=\bar L\, \displaystyle{\frac{\p U(L, \bar L)}{\p \bar L}}
\eeq

Notice that (\ref{canon}) gives an equivalent description
of the canonical pair from the Ta\-ka\-sa\-ki-\-Ta\-ke\-be
theorem.
Namely, (\ref{canon}) tells us that
$U(L, \bar L)$ is {\it generating function of the canonical
transformation}
$L, M \longrightarrow \bar L^{-1}, \bar M$ and, therefore,
implicitly determines the corresponding canonical $f,g$-pair.
Indeed, let in general
$M=F(L, \bar L)$,
$\bar M=G(L, \bar L)$,
and find the Poisson bracket between $\bar L^{-1}$ and $\bar M$.
For brevity we write $\bar L = \varphi (L,M)$, where $\varphi$
is a function implicitly determined
by the equality $M=F(L, \bar L)$. We have:
$$
\frac{\p \bar L}{\p L}
\frac{\p \bar M}{\p M} -
\frac{\p \bar L}{\p M}
\frac{\p \bar M}{\p L} =
\frac{\p \varphi}{\p L}
\frac{\p G}{\p \bar L}
\frac{\p \varphi}{\p M} -
\frac{\p \varphi}{\p M}
\left (
\frac{\p G}{\p \bar L}+
\frac{\p G}{\p \bar L}
\frac{\p \varphi}{\p L}\right ) =
-\, \frac{\p \varphi}{\p M}
\frac{\p G}{\p \bar L}
$$
Therefore,
$$
\{ \bar L^{-1}, \, \bar M \} = L\bar L^{-2}
\frac{\p \varphi}{\p M} \frac{\p G}{\p \bar L}
= \frac{L}{\bar L^2}\,
\frac{\p G/\p L}{\p F/ \p \bar L}
$$
where the last equality follows from
$\p \varphi / \p M =(\p \bar L/\p M)_{L=const} =(\p F/\p \bar L)^{-1}$.
Now, on taking $F, G$ from (\ref{canon}), we get
$\{ \bar L^{-1}, \, \bar M \} =\bar L^{-1}$,
i.e. the transformation (\ref{canon}) is indeed canonical.

The Poisson bracket
$\{ L, \, \p_{L} U(L, \bar L) \} =1$
follows from
$\{ L, \, M \} =L$. This can be rewritten as
\beq\label{string2}
\{ L, \, \bar L \} =\frac{1}{\sigma (L, \bar L)}
\eeq
which is one of the forms of the so-called string equation.

\section{Conclusion}

We have demonstrated that some classical problems of
complex analysis
have an integrable structure which is the same as
the one arising in the dispersionless limit
of soliton equations.
Given a simply connected
domain in the complex plane
and an auxiliary real-analytic
function $\sigma (z,\bar z)$ defined
all over the complex plane, we have constructed
a solution to the dispersionless 2D Toda chain
hierarchy. The independent flows of the hierarchy
are moments of the domain, or multipole
coefficients, defined with respect to
the function $\sigma$ which in the context
of potential theory plays the role of density
of background charge. The set of dependent variables
is then the complimentary set of multipole coefficients.
The energy of the background charge in the domain,
regarded as a function of the moments,
is (logarithm of) the $\tau$-function of the
hierarchy. This function obeys the dispersionless
Hirota equations.
Taking derivatives of this function
with respect to the moments, one obtains
formal solutions to the conformal mapping
problem, the Dirichlet boundary problem, and
to the 2D inverse potential problem.
An interesting question, to be discussed elsewhere,
is an inverse problem connected with the scheme
outlined above: given a solution to the dispersionless
Hirota equations,
associate with it a boundary or conformal mapping problem,
and to restore the function $\sigma$.
Such an interpretation may
be helpful for better understanding dispersionless
integrable hierarchies.

\section*{Acknowledgements}

The author is grateful to the organizers of the
Protvino conference ``Classical and quantum
integrable systenms'' for the opportunity to present
these results, and to I.Kri\-che\-ver, A.Mar\-sha\-kov,
M.Mi\-ne\-ev-\-Wein\-stein
and P.Wi\-eg\-mann for useful discussions.
This work was supported in part by CRDF grant RP1-2102,
by grant INTAS-99-0590 and by RFBR grant 00-02-16477.

\section*{Appendix}

In the appendix we give a direct proof of the formulas
(\ref{25}) for first derivatives of the function $F$
with respect to $t_k$. We use the notation of Sec.\,2.

Consider the variation $\delta_k F$ ($k>0$) of the function
$$
F=\frac{1}{2\pi}\int_{D_{+}}\tilde \Phi (z,\bar z)
\sigma (z, \bar z)d^2z
$$
under a small change $t_k \to t_k +\delta t_k$:
$F\to F+\delta_k F$. All other
moments $t_j$ with $j\neq k$, $j\geq 0$, are kept constant.
We have:
$$
 \delta_k F = \frac{1}{2\pi}\int_{D_{+}}
\delta_k \tilde \Phi^{+} (z,\bar z)\sigma (z, \bar z) d^2z +
\frac{1}{2\pi}\int_{\delta_k D_{+}}
\tilde \Phi (z,\bar z) \sigma (z, \bar z)d^2z
$$
where $\delta_k D_{+}$ is the variation of the domain
due to the change of $t_k$'s. (To be more precise,
$\int_{\delta_k D_{+}}$ stands for the difference
$\int_{ D_{+}(t_k +\delta t_k)}- \int_{D_{+}(t_k)}$.)
By $I_1$ and $I_2$ we denote the first and the second
terms of this formula. $I_1$ is easy to calculate. Since
$\delta_k \tilde \Phi^{+} (z,\bar z)=
z^k \delta t_k +\bar z^k \delta \bar t_k$,
the first term is
$$
I_1 = \frac{1}{2\pi}\int_{D_{+}}
\delta_k \tilde \Phi^{+} (z,\bar z) \sigma (z, \bar z)d^2z =
\frac{1}{2}(v_k \delta t_k +\bar v_k \delta \bar t_k)
$$

Let us turn to the second one,
$I_2 = \frac{1}{2\pi}\int_{\delta_k D_{+}}
\tilde \Phi (z,\bar z) \sigma (z, \bar z)d^2z$.
Since the integral
goes over a small neighbourhood of the curve
$\gamma =\p D_{+}$, we can substitute $\tilde \Phi$  by
$\tilde \Phi^{-}$,
which, by virtue of (\ref{cont}),
coincides in this neighbourhood
with $\tilde \Phi^{+}$ up to second order terms.
Recall that
$\tilde \Phi^{-} (z,\bar z)= v_0 +2\,{\cal R}e \,\omega (z)$,
where the function
$$
\omega (z)=\sum_{k\geq 1} \frac{v_k}{k}z^{-k}
$$
is holomorphic in $D_{-}$ and regular at infinity.
Under our assumptions,
all singularities of $\omega$ are in a domain
$B\subset D_{+}$ such that
$R= D_{+}\backslash B$ is a strip-like domain.
In other words,
$\omega (z)$ can be analytically extended
across $\gamma$ to the strip-like domain $R$.
Moreover, recall that we are working
in the space of analytic curves. It can be proved that
for curves that belong to
a small but finite neighbourhood
of the point in this space that corresponds to a given curve,
the domain $B$ can be chosen to be independent of $t_k$.

We can write, in the first order in $\delta t_k$:
$$
I_2 =
\frac{v_0}{2\pi}\int_{\delta_k D_{+}}
\sigma (z, \bar z) d^2z +
\frac{1}{\pi} {\cal R}e \int_{\delta_k R}
\omega (z) \sigma (z, \bar z)d^2z
$$
The first integral in the right hand side
is the variation of $t_0$, and, therefore, equals zero
as soon as $t_0$ is kept constant.
So, we are left with
the second one,
$$
\begin{array}{lll}
I_2 &=&
\displaystyle{
\frac{1}{\pi} {\cal R}e \int_{D_{+}(t_k +\delta t_k )
\backslash B}\!
\omega (z) \sigma (z,\bar z)d^2z
-\frac{1}{\pi} {\cal R}e \int_{D_{+}(t_k )
\backslash B}\!
\omega (z) \sigma (z, \bar z)d^2z }\, =
\\ && \\
&=&
\displaystyle{
{\cal R}e \oint_{\p R(t_k +\delta t_k )}\!
\frac{dz}{2\pi i}\,
\omega (z) \p_z U(z,\bar z)
- {\cal R}e \oint_{\p R(t_k )}\!
\frac{dz}{2\pi i} \,
\omega (z) \p_z U(z, \bar z)}
\end{array}
$$
The second line follows from the Stokes formula.
The boundary of $R$ consists of two curves: inner
and outer ones. The inner boundary is the same
in both integrals, so its contribution cancel.
The outer boundary is $\gamma (t_k +\delta t_k )$
in the first integral and $\gamma (t_k)$ in the second one.
Therefore,
using the definition of the generalized
Schwarz function (\ref{sch1}),
we can write:
$$
I_2 =
{\cal R}e \oint_{\gamma (t_k +\delta t_k )}\!
\frac{dz}{2\pi i}\,
\omega (z) S(t_k \!+ \! \delta t_k , z)
- {\cal R}e \oint_{\gamma (t_k )}\!
\frac{dz}{2\pi i}\,
\omega (z) S(t_k , z)
$$
Since $S(z)$ is holomorphic in $R$, the contour
can be taken to be $\gamma =\gamma (t_k)$ in both integrals.
Thus we obtain, by virtue of (\ref{sch2}):
$$
I_2 = {\cal R}e \left ( \delta t_k
\oint_{\gamma}\! \frac{dz}{2\pi i} \,
\omega (z) \frac{\p S(z)}{\p t_k} \right )=
\frac{1}{2}(v_k \delta t_k + \bar v_k \delta \bar t_k )
$$
Summing up the contributions $I_1$ and $I_2$, we get
$\delta_k F= v_k \delta t_k + \bar v_k \delta \bar t_k$,
which is the desired result for $k\geq 1$.
The calculation of the $t_0$-derivative can be done
in a similar way.

\end{document}